\documentclass{amsart}
 \usepackage{amssymb,amsmath,amsfonts,epsfig,latexsym}
 
 \newtheorem{theorem}{Theorem}[section]

\newtheorem{proposition}[theorem]{Proposition}
\newtheorem{lemma}[theorem]{Lemma}
\newtheorem{corollary}[theorem]{Corollary}

\newtheorem*{BF}{Braun's Formula}

\theoremstyle{definition}
\newtheorem{remark}[theorem]{Remark}
\newtheorem{example}[theorem]{Example}

\def\Q{\mathbb{Q}} 
\def\R{\mathbb{R}}
\def\T{\mathcal{T}}
\def\Z{\mathbb{Z}}

\def\sing{\mathrm{sing}}

\def\<{\langle}
\def\>{\rangle}
 
\DeclareMathOperator{\BBox}{Box}

\DeclareMathOperator{\conv}{conv}
\DeclareMathOperator{\Ehr}{Ehr}
\DeclareMathOperator{\Hom}{Hom}

\DeclareMathOperator{\lk}{lk}

\DeclareMathOperator{\Star}{Star}

\DeclareMathOperator{\vol}{vol}

\begin{document}

\title{Ehrhart series and lattice triangulations}           
    
\author{Sam Payne}
\address{Stanford University, Mathematics, Bldg. 380, 450 Serra Mall, Stanford, CA 94305}
\email{spayne@stanford.edu}
\thanks{Supported by the Clay Mathematics Institute.}

\begin{abstract}
We express the generating function for lattice points in a rational polyhedral cone with a simplicial subdivision in terms of multivariate analogues of the $h$-polynomials of the subdivision and ``local contributions" of the links of its nonunimodular faces.  We also compute new examples of nonunimodal $h^*\!$-vectors of reflexive polytopes.
\end{abstract}
 
\maketitle

\section{Introduction}

Let $N$ be a lattice and let $\sigma$ be a strongly convex rational polyhedral cone in $N_\R = N \otimes_\Z \R$.  The generating function for lattice points in $\sigma$,
\[
G_\sigma = \sum_{v \in (\sigma \cap N)} x^v,
\]
is a rational function, in the quotient field $\Q(N)$ of the multivariate Laurent polynomial ring $\Z[N]$.  For instance, if $\sigma$ is unimodular, spanned by a subset $\{v_1, \ldots, v_r\}$ of a basis for $N$, then $G_\sigma$ is equal to $1/(1-x^{v_1}) \cdots (1-x^{v_r})$.

Suppose $\Delta$ is a rational simplicial subdivision of $\sigma$, with $v_1, \ldots, v_s$ the primitive generators of the rays of $\Delta$.  We define a multivariate analogue $H_\Delta$ of the $h$-polynomial $h_\Delta(t) = \sum_{\tau \in \Delta} t^{\dim \tau}\cdot (1-t)^{\dim \sigma - \dim \tau}$, by
\[
H_\Delta = \sum_{\tau \in \Delta} \bigg( \prod_{v_i \in \tau} x^{v_i} \cdot \prod_{v_j \not \in \tau} (1-x^{v_j}) \bigg).
\]
Every point in $\sigma$ is in the relative interior of a unique cone in $\Delta$ and can be written uniquely as a nonnegative integer linear combination of the primitive generators of the rays of that cone plus a fractional part.  The generating function for lattice points in the relative interior of $\tau$ that have no fractional part is $\prod_{v_i \in \tau} \, x^{v_i} / (1-x^{v_i})$, so the generating function for lattice points in $\sigma$ with no fractional part (with respect to the subdivision $\Delta$) is $H_\Delta / (1-x^{v_1}) \cdots (1-x^{v_s})$.  In particular, if every cone in $\Delta$ is unimodular, spanned by part of a basis for the lattice, then every lattice point has no fractional part, so this gives $G_\sigma$.  Otherwise, the remaining lattice points with nonzero fractional part, which necessarily lie in the nonunimodular cones of $\Delta$, may be accounted for as follows.

Say that a cone is singular if it is not unimodular, and let $\Delta^\sing$ be the set of singular cones in $\Delta$.  After possibly renumbering, say $\tau \in \Delta^\sing$ is spanned by $\{ v_1, \ldots, v_r \}$.  Let $\BBox(\tau)$ be the open parallelipiped
\[
\BBox(\tau) = \{ a_1 v_1 + \cdots + a_r v_r \ | \ 0 < a_i < 1 \},
\]
and let $B_\tau$ be the generating function for lattice points in $\BBox(\tau)$,
\[
B_\tau\,  = \sum_{v \, \in \, \BBox(\tau) \cap N} x^v.
\]
We write $\lk \tau$ for the link of $\tau$ in $\Delta$.  In other words, $\lk \tau$ is the union of the cones $\gamma$ in $\Delta$ such that $\gamma \cap \tau = 0$ and $\gamma + \tau$ is a cone in $\Delta$.  We define a multivariate analogue $H_{\lk \tau}$ of the $h$-polynomial of $\lk \tau$ by
\[
H_{\lk \tau} = \sum_{\gamma \in \lk \tau} \bigg( \prod_{v_i \in \gamma} x^{v_i} \cdot \prod_{v_j \in (\lk \tau \smallsetminus \gamma) } (1-x^{v_j}) \bigg).
\]

\vspace{5 pt}

\begin{theorem} \label{main}
Let $\Delta$ be a rational simplicial subdivision of a strongly convex rational polyhedral cone $\sigma$, and let $v_1, \ldots, v_s$ be the primitive generators of the rays of $\Delta$.  Then
\[
(1-x^{v_1}) \cdots (1-x^{v_s}) \cdot G_\sigma \ = \ H_\Delta + \sum_{\ \tau \, \in \, \Delta^{\sing}} \bigg( B_\tau \cdot H_{\lk \tau} \cdot \prod_{v_i \not \in \Star \tau} (1-x^{v_i}) \bigg).
\]
\end{theorem}

\noindent Here $\Star \tau$ is the union of the maximal cones in $\Delta$ that contain $\tau$.

Recent work of Athanasiadis \cite{Athanasiadis04, Athanasiadis05}, Bruns and R\"omer \cite{BrunsRomer07}, and Ohsugi and Hibi \cite{OhsugiHibi06} has highlighted the usefulness of considering the effects of a ``special simplex" that is contained in all of the maximal faces of a triangulation.   Our next result, which is inspired by their work, is a generalization of Theorem \ref{main} that takes into account the effect of a special cone that is contained in all of the maximal cones of the subdivision.  For a cone over a triangulated polytope, the notion of special cone that we consider is slightly more general than a cone over a special simplex in the sense of \cite{Athanasiadis04}; a special cone is a cone over a special simplex if and only if it is not contained in the boundary of the cone over the polytope.

Suppose, as above, that $\tau$ is a cone in $\Delta$ spanned by $v_1, \ldots, v_r$, and let $\lambda$ be a face of $\tau$.  After possibly renumbering, we may assume that $\lambda$ is spanned by $v_1, \ldots, v_q,$ for some $q \leq r$.  We define a partially open parallelipiped $\BBox(\tau, \lambda)$, which we think of as the $\BBox$ of $\tau$ relative to $\lambda$, as
\[
\BBox(\tau, \lambda) = \{ a_1 v_1 + \cdots + a_r v_r \ | \ 0 \leq a_i < 1 \mbox{ for all } i, \mbox{ and } a_i \neq 0 \mbox{ for } i > q \}.
\]
Let the polynomial $B_{\tau, \lambda}$ be the relative analogue of $B_\tau$,
\[
B_{\tau, \lambda} = \sum_{v \, \in \, \BBox(\tau, \lambda) \cap N} x^v.
\]

\begin{theorem} \label{special cone}
Let $\Delta$ be a rational simplicial subdivision of a strongly convex rational polyhedral cone $\sigma$, and let $v_1, \ldots, v_s$ be the primitive generators of the rays of $\Delta$.  Suppose $\lambda \in \Delta$ is contained in every maximal cone of $\Delta$.  Then
\[
(1-x^{v_1}) \cdots (1-x^{v_s}) \cdot G_\sigma \ = \ H_{\lk \lambda} + \sum_{\tau \succeq \lambda} \bigg( B_{\tau, \lambda} \cdot H_{\lk \tau} \cdot  \prod_{v_j \not \in \Star \tau}(1-x^{v_j}) \bigg).
\]
\end{theorem}

\noindent Since the zero cone is contained in every cone of $\Delta$ and $B_{\tau,0} = B_\tau$ for all $\tau$, Theorem~\ref{main} is the special case of Theorem~\ref{special cone} where $\lambda = 0$.

\begin{remark}
There are many ways of computing the rational function $G_\sigma$ and its specializations, some of which are algorithmically efficient.  The excellent survey articles \cite{BarvinokPommersheim99} and \cite{DeLoera05} may serve as introductions to the extensive literature on this topic.  Efficient algorithms have been implemented in the computer program LattE \cite{DeLoera04}.  Theorems~\ref{main} and \ref{special cone}, and their specializations to lattice polytopes (Corollaries~\ref{Batyrev and Dais} and \ref{special simplex}), seem to be useful in cases where it is especially easy to give subdivisions that are close to unimodular, and in studying families of such examples in which the contributions of the singular cones can be easily understood.  See, for instance, the examples in Section \ref{examples}, which were computed by hand.
\end{remark}

Theorems \ref{main} and \ref{special cone} specialize to give formulas for Ehrhart series of lattice polytopes.  See Section~\ref{specialization} for details.  We use these specializations to construct examples of reflexive polytopes with interesting $h^*\!$-vectors, the results of which are summarized as follows.  Suppose $N'$ is a lattice and $P$ is a $d$-dimensional lattice polytope in $N'_\R$.  Recall that the Ehrhart series of $P$ is 
\[
\Ehr_P(t) = 1 + \sum_{m \geq 1} \#\{ mP \cap N'\} \cdot t^m,
\]
and that $(1-t)^{d+1} \Ehr_P(t) = h^*_0 + \cdots + h^*_d t^d$ for some integers $h^*_i$.  See \cite[Chapters 3 and 4]{BeckRobins07} for these and other basic facts about Ehrhart series.  We say that $h^*(P) = (h^*_0, \ldots, h^*_d)$ is the $h^*\!$-vector of $P$, and write $h^*_P(t)$ for the polynomial $\sum h^*_i t^i$.\footnote{There is unfortunately no standard notation for the numerator of the Ehrhart series of a lattice polytope and its coefficients, despite the extensive literature on the topic.  The relation between this polynomial and the Ehrhart polynomial of a lattice polytope is analogous to the relation of the $h$-polynomial of a simplicial polytope with the $f$-polynomial; some like to call its vector of coefficients the ``Ehrhart $h$-vector" \cite{BrunsRomer07}.  Others simply denote it by $h$ \cite{Athanasiadis05, BeckRobins07, OhsugiHibi06} or $\delta$ \cite{Hibi90, Hibi91, Hibi92, Hibi94}.  Here we use the notation $h^*$, following Stanley \cite{Stanley93} and Athanasiadis \cite{Athanasiadis04}, to emphasize the analogy with $h$-polynomials of simplicial polytopes while avoiding any possible ambiguity.}

Recall that $P$ is said to be reflexive if it contains $0$ in its interior and its dual polytope has vertices in the dual lattice.  Suppose that $P$ is reflexive.  Then $h^*_i = h^*_{d-i}$ for all $i$ \cite{Hibi90}.  Furthermore, if $\partial P$ has a regular unimodular triangulation, given by intersecting $\partial P$ with the domains of linearity of a convex piecewise linear function on $N'_\R$, then $h^*(P)$ is equal to the $h$-vector of this triangulation, which is combinatorially equivalent to the boundary complex of a simplicial polytope.  In particular, if $\partial P$ has a regular unimodular triangulation then $h^*(P)$ is unimodal, in the sense that $h^*_0 \leq h^*_1 \leq \cdots \leq h^*_{[d/2]}$, and furthermore the vector of successive differences 
\[
g^*(P) = (h^*_0, h^*_1 - h^*_0, \ldots, h^*_{[d/2]} - h^*_{[d/2] - 1})
\]
is a Macaulay vector, i.e.\ the Hilbert sequence of a graded algebra generated in degree one.  For an arbitrary reflexive polytope, Hibi showed that $h^*_0 \leq h^*_1 \leq h^*_i$ for $2 \leq i <d$, so $h^*(P)$ is unimodal if $d \leq 5$ and $g^*(P)$ is a Macaulay vector if $d \leq 3$ \cite{Hibi91}.

First examples of reflexive polytopes with nonunimodal $h^*\!$-vectors were given in \cite{MustataPayne05} in even dimensions $d \geq 6$.  These examples were nonsimplicial, the depth of the ``valleys" in $h^*(P)$ was never more than two, and the construction did not yield nonunimodal examples in odd dimensions.  Also, it remained unclear, in the cases where $h^*(P)$ is unimodal, whether $g^*(P)$ is necessarily a Macaulay vector.

\begin{theorem} \label{not unimodal}
For every $d \geq 6$, there exists a $d$-dimensional reflexive simplex $P$ such that $h^*(P)$ is not unimodal.
\end{theorem}

\begin{theorem} \label{not macaulay}
For every $d \geq 4$, there exists a $d$-dimensional reflexive simplex $P$ such that $h^*(P)$ is unimodal, but $g^*(P)$ is not a Macaulay vector.
\end{theorem}

\begin{theorem} \label{big valleys}
For any positive integers $m$ and $n$, there exists a reflexive polytope $P$ and indices $i_1 < j_1< i_2< j_2 < \cdots < i_m < j_m < i_{m+1}$ such that
\[
h^*_{i_\ell} - h^*_{j_\ell} \geq n \mbox{ \ and \ } h^*_{i_{\ell+1}} - h^*_{j_\ell} \geq n,
\]
for $1 \leq \ell \leq m$.  Furthermore, $P$ can be chosen so that $\dim(P) = O(m \log \log n)$.
\end{theorem}

\noindent In other words, for any positive integers $m$ and $n$, there exists a reflexive polytope $P$ of dimension $O(m \log \log n)$ such that $h^*(P)$ has at least $m$ valleys of depth at least $n$.

\begin{remark}
None of the examples of reflexive polytopes with nonunimodal $h^*\!$-vectors constructed here are normal, in the sense where a lattice polytope $P$ is normal if every lattice point in $mP$ is a sum of $m$ lattice points in $P$, for all positive integers $m$.  For normal reflexive polytopes $P$, the questions of whether $h^*(P)$ is unimodal and whether $g^*(P)$ is a Macaulay vector remain open and interesting \cite{OhsugiHibi06}.
\end{remark}

We conclude the introduction with an example illustrating Theorem~\ref{special cone}.

\begin{example}
Suppose $N = \Z^3$, $v = (0,0,1)$ and
\[
\begin{array}{llll}
 v_1 = (1,0,1), & v_2 = (0,1,1), & v_3 = (0,-1,1), & v_4 = (-1,0,1).
\end{array}
\]
Let $\sigma$ be the cone spanned by $v_1, \, v_2,\,  v_3,$ and $v_4$, with $\Delta$ the simplicial subdivision of $\sigma$ whose maximal cones are
\[
\sigma_1 = \< v_1, v_2, v_3 \> \mbox{  \ and \ } \sigma_2 = \< v_2, v_3, v_4 \>.
\]
Then $\sigma_1$ and $\sigma_2$ are not unimodular; the lattice points with no fractional part are exactly those $(a,b,c)$ such that $a + b + c$ is even.  The remaining lattice points may be written uniquely as $v_0$ plus a lattice point with no fractional part.  We compute the generating function $G_\sigma$ as follows, using Theorem~\ref{special cone} with $\lambda = \<v_2 , v_3\>$.

There are exactly three cones in $\lk \lambda$: 0, $\<v_1\>$, and $\< v_4\>$.  Therefore
\begin{eqnarray*}
H_{\lk \lambda} &=& (1-x^{v_1}) (1-x^{v_4}) + x^{v_1} (1-x^{v_4}) + x^{v_4}(1-x^{v_1}). \\
  &=& 1 - x^{v_1}x^{v_4}.
  \end{eqnarray*}
    Since $\BBox_{\lambda, \lambda}$ contains a unique lattice point $v$ and $\BBox_{\sigma_1, \lambda}$ and $\BBox_{\sigma_2, \lambda}$ contain no lattice points, it follows that
\[
G_\sigma = \frac{(1 + x^{v})(1-x^{v_1}x^{v_4})}{(1-x^{v_1})(1-x^{v_2})(1-x^{v_3})(1-x^{v_4})}.
\]
\end{example}

\vspace{5 pt}

\noindent \textbf{Acknowledgments.}  This note was written in response to questions raised at the 2006 Snowbird conference on Integer Points in Polyhedra.  I am grateful to the organizers of the conference for the opportunity to participate and thank A.~Barvinok, M.~Beck, B.~Braun, J.~De Loera, C.~Haase, T.~Hibi, and S.~Ho\c{s}ten for stimulating conversations.  I also thank P.~McMullen and the referee for helpful comments and corrections.

\section{Proof of Theorems 1.1 and 1.2}

As observed in the introduction, Theorem~\ref{main} is a special case of Theorem~\ref{special cone}.  We will begin by showing the converse, that Theorem~\ref{special cone} is a consequence of Theorem~\ref{main}, using the following lemma.  The lemma is a multivariate analogue of the familiar fact that, for any simplcial complex $\Delta'$, the $h$-polynomial of the join of $\Delta'$ with a simplex is equal to the $h$-polynomial of $\Delta'$.

\begin{lemma}  \label{H identity}
If $\lambda$ is a face of every maximal cone of $\Delta$ then $H_{\Delta} = H_{\lk \lambda}$ and
\[
H_{\lk \gamma} = H_{\lk (\gamma + \lambda)},
\]
for every $\gamma \in \Delta$.
\end{lemma}

\begin{proof}
We show that $H_\Delta = H_{\lk \lambda}$.  The proof of the second claim is similar.  Every maximal face of $\Delta$ contains $\lambda$ if and only if $\Delta$ is the join of $\lambda$ with $\lk \lambda$.  Since $\lambda$ is the join of its rays, it will suffice to consider the case where $\lambda$ is one-dimensional, with primitive generator $v_1$.  In this case, the required identity may be seen by regrouping the terms in the summation defining $H_\Delta$ as
\[
H_\Delta = \sum_{v_1 \not \in \gamma} \bigg( \prod_{v_i \in \gamma} x^{v_i} \cdot \prod_{v_j \not \in \gamma} (1-x^{v_j}) \bigg) +  \sum_{v_1 \in \tau} \bigg( \prod_{v_k \in \tau} x^{v_k} \cdot \prod_{v_\ell \not \in \tau} (1-x^{v_\ell}) \bigg). 
\]
Since the cones not containing $v_1$ are exactly the $\gamma \in \lk \lambda$, and since the cones containing $v_1$ are exactly those $\tau = (\gamma + \lambda)$ for $\gamma \in \lk \lambda$, the above equation gives
\[
H_\Delta \ = \ (1-x^{v_1}) \cdot H_{\lk \lambda} \ + \ x^{v_1} \cdot H_{\lk \lambda},
\]
so $H_\Delta = H_{\lk \lambda}$, as required.
\end{proof}

\noindent  Suppose $\lambda$ is contained in every maximal cone.  Then, for any cone $\tau$ containing $\lambda$,
\[
B_{\tau, \lambda} = \sum_{(\gamma + \lambda) = \tau} B_\gamma.
\]
Futhermore, for each $\gamma$ such that $(\gamma + \lambda) = \tau$, we have $\Star \gamma = \Star \tau$, and $H_{\lk \gamma} = H_{\lk \tau}$, by Lemma~\ref{H identity}.  Therefore, Theorem~\ref{special cone} follows from Theorem~\ref{main}, which we now prove.

\begin{proof}[Proof of Theorem~\ref{main}]
Let $v$ be a lattice point in $\sigma$.  Then $v$ is contained in the relative interior of a unique cone $\gamma \in \Delta$.  If $v_1, \ldots, v_r$ are the primitive generators of the rays of $\gamma$, then $v$ can be written uniquely as 
\[
v = a_1v_1 + \cdots + a_r v_r + \{v \},
\]
where each $a_i$ is a nonnegative integer and $\{ v \}$, which we call the fractional part of $v$, is either zero or lies in $\BBox(\tau)$ for some unique singular cone $\tau \preceq \gamma$.  Since $v$ lies in the relative interior of $\gamma$, $a_i$ must be strictly positive for each $v_i$ in $(\gamma \smallsetminus \tau)$.  Conversely, if $v'$ is a lattice point in $\BBox(\tau)$, and $v = a_1 v_1 + \cdots + a_r v_r + v'$, where the $a_i$ are nonnegative integers that are strictly positive for $v_i \in (\gamma \smallsetminus \tau)$, then $v \in \gamma^\circ$ and $\{v \} = v'$.  Therefore,
\[
\frac{x^{v'} \cdot H_{\lk \tau}}{\prod_{v_i \in \Star{\tau}} (1-x^{v_i})}
\]
is the generating function for lattice points $v \in \sigma$ such that $\{ v \} = v'$. Then $B_\tau \cdot H_{\lk \tau} / \prod_{v_i \in \Star \tau} (1-x^{v_i})$ is the generating function for lattice points in $\sigma$ whose fractional part is in $\BBox(\tau)$, and the theorem follows.
\end{proof}

\section{Specialization to lattice triangulations of polytopes} \label{specialization}

Let $N'$ be a lattice, and let $P$ be a $d$-dimensional lattice polytope in $N'_\R$.  Suppose $N = N' \times \Z$, with $\sigma$ the cone over $P \times \{1\}$ in $N_\R$, and let $u : N \rightarrow \Z$ be the projection to the second factor.  Since
\[
\#\{mP \cap N'\} =   \# \{ v \in (\sigma \cap N) \, | \, u(v) = m \}
\]
for all positive integers $m$, the specialization
\[
\varphi : \Q(N) \rightarrow \Q(t), \ \ x^v \mapsto t^{u(v)}
\]
maps $G_\sigma$ to $\Ehr_P(t)$.

Suppose $\T$ is a lattice triangulation of $P$, and let $\Delta$ be the subdivision of $\sigma$ consisting of the cones $\tau_F$ over $F \times \{1\}$ for all faces $F \in \T$.  Then each of the primitive generators $v_1, \ldots, v_s$ of the rays of $\Delta$ is a lattice point in $P \times \{1 \}$, so $u(v_i) = 1$ for all $i$.  It follows that
\[
\varphi (H_\Delta) = (1-t)^{(s-d-1)} \cdot h_\T(t),
\]
where $h_\T(t)$ is the $h$-polynomial of the simplicial complex $\T$.  Similarly,
\[
\varphi (H_{\lk \tau_F}) = (1-t)^{(s' - d' -1)} \cdot h_{\lk F} (t),
\]
where $s'$ is the number of vertices in $\lk F$, and $d' = \dim \lk F$.  Let $\T^\sing$ denote the set of nonunimodular simplices of $\T$.  For $F \in \T^\sing$, let $B_F(t) = \varphi (B_{\tau_F})$, so
\[
B_F(t) = \sum_{v \in \BBox(\tau_F) \cap N} t^{u(v)}.
\]
The following specialization of Theorem~\ref{main} is due to Betke and McMullen \cite[Theorem~1]{BetkeMcMullen85} and was rediscovered by Batyrev and Dais in the context of stringy algebraic geometry  \cite[Theorem~6.10]{BatyrevDais96}.

\begin{corollary} \label{Batyrev and Dais}  Let $P$ be a lattice polytope, and let $\T$ be a lattice triangulation of $P$.  Then
\[
h^*_P(t) = h_\T(t) + \sum_{F \in \T^\sing} B_F(t) \cdot h_{\lk F} (t).
\]
\end{corollary}

\noindent Note that $B_F(t)$ and $h_{\lk F}(t)$ have nonnegative integer coefficients, so the theorem of Stanley that $h^*_P(t)$ has nonnegative integer coefficients  \cite[Theorem~2.1]{Stanley80} follows immediately.  Another interesting proof of this nonnegativity, using ``irrational decompositions," recently appeared in work of Beck and Sottile \cite{BeckSottile06}.  

Furthermore, $h_{\lk F}(t)$ is always nonzero, and $B_F(t)$ must be nonzero for some $F$ if $\T$ is not unimodular.  Therefore, from Corollary \ref{Batyrev and Dais}, we deduce the following.

\begin{corollary}  \cite[Theorem~2]{BetkeMcMullen85}
Let $P$ be a lattice polytope, and let $\T$ be a lattice triangulation of $P$.  Then
\[
h^*_i(P) \geq h_i(\T)
\]
for all $i$.  Furthermore, equality holds for all $i$ if and only if $\T$ is unimodular.
\end{corollary}

If the triangulation $\T$ contains a special simplex $F'$ that is a face of every maximal simplex of $\T$, then we can take this into account using Theorem~\ref{special cone}.  For $F \geq F'$, let $B_{F,F'}(t) = \varphi (B_{\tau_F, \tau_{F'}})$.

\begin{corollary} \label{special simplex}
Let $P$ be a lattice polytope, and let $\T$ be a triangulation of $P$ with a special simplex $F'$.  Then
\[
h^*_P(t) = h_{\lk F'}(t) + \sum_{F \geq F' } B_{F,F'}(t) \cdot h_{\lk F}(t).
\]
\end{corollary}

\noindent Recall that a lattice polytope $P$ is called reflexive if it contains $0$ in its interior and the polar dual polytope of $P$ has vertices in the dual lattice of $N'$.  In the special case where $P$ is reflexive and $0 = F'$ is a special simplex of the triangulation $\T$, then $\T$ is the join of $\{0\}$ with a triangulation of the boundary of $P$, and we recover \cite[Theorem~1.3]{MustataPayne05}.

\begin{remark}
In \cite{BrunsRomer07}, Bruns and R\"omer use techniques from commutative algebra to show that if $mP$ is a translate of a reflexive polytope for some positive integer $m$, and if $P$ has any regular unimodular triangulation, then $P$ has a regular unimodular triangulation with an $(m-1)$-dimensional special simplex $F'$ such that $\lk F'$ is combinatorially equivalent to the boundary complex of a simplicial polytope.  From this they deduce that $h^*(P) = h(\lk F')$ is unimodal and that $g^*(P)$ is a Macaulay vector.
\end{remark}

\begin{remark}
These specializations to lattice polytopes are closely related to stringy invariants of toric varieties.  See \cite{BatyrevDais96} and \cite{MustataPayne05} for details on this connection.  Special simplices have also appeared in the stringy geometry literature, where they have been called the ``core" of a triangulation \cite{Stienstra98}.
\end{remark}

\section{Examples of reflexive polytopes with nonunimodal $h^*$-vectors} \label{examples}

Let $N'$ be a lattice and let $M' = \Hom(N',\Z)$ be its dual lattice.  Let $P$ be a $d$-dimensional lattice polytope in $N'_\R$.  Recall that $P$ is reflexive if and only if it contains $0$ in its interior and the dual polytope
\[
P^\circ = \{ u \in M'_\R \ | \ \<u, v \> \geq -1 \mbox{ for all } v \in P \}
\]
has vertices in $M'$.

\begin{proposition} \label{simplices}
Let $a_1, \ldots, a_d,$ and $b$ be positive integers, and let
\[
f = \frac{1}{b}(a_1, \ldots, a_d).
\]
Then the simplex $P = \conv \{ e_1, \ldots, e_d, -f\}$ in $\R^d$ is reflexive with respect to the lattice $N' = \Z^d + \Z \cdot f$ if and only if $a_1 + \cdots + a_d = bc$ for some integer $c$ and each $a_i$ divides $b(c + 1)$.
\end{proposition}

\begin{proof}
Let $u = -e_1^* - \cdots - e_d^*$.  Then $P^\circ$ is the simplex in $(\R^d)^*$ with vertices $u$ and
\[
u_i = u + \frac{a_1 + \cdots + a_d + b}{a_i} e_i^*,
\]
for $1 \leq i \leq d$.  Now $u$ lies in $M'$ if and only if $\<u, f\> = (a_1 + \cdots + a_d)/b$ is some integer $c$.  If this is the case, then $u_i$ lies in $M'$ if and only if $\<u_i, e_i \> = b(c+1) / a_i - 1$ is an integer.
\end{proof}

In the following examples, we use Proposition~\ref{simplices} to construct a reflexive simplex $P$, and then apply Corollary~\ref{special simplex} with respect to the triangulation $\T$ obtained by taking the join of $\{0\}$ with the boundary of $P$, and the special simplex $F' = 0$, to compute $h^*(P)$.  The link of a $(d-r)$-dimensional face $F \in \T$ that contains $0$ is the boundary complex of an $r$ simplex, so $h^*_{\lk F} (t) = 1 + t + \cdots + t^r$.  We write $\widetilde{v} = (v,1)$ for the lift of a point $v \in \R^d$ to height one in $\R^d \times \R$.

\begin{example} \label{reflexive simplex}
Suppose $d = bk + r$ for some positive integers $b$ and $k$ and some nonnegative integer $r$.  Let $f$ be the following vector in $\R^d$,
\[
f = \frac{1}{b}(\underbrace{1, \ldots, 1}_{bk}, \underbrace{b, \ldots, b}_r ).
\]
By Proposition~\ref{simplices}, $P = \conv\{ e_1, \ldots, e_d, -f \}$ is reflexive with respect to the lattice $N'= \Z^d + \Z \cdot f$. The only nonunimodular face $F$ of $\T$ such that $\BBox(\tau_F, 0)$ contains lattice points is $F = \conv \{0, e_1, \ldots, e_{bk}\}$.  Let $f' = f - e_{bk + 1} - \cdots - e_d$.  Then $\BBox(F,0) \cap (N' \times \Z)$ consists of the lattice points
\[
\frac{1}{b}(\widetilde e_1 + \cdots + \widetilde e_{bk}), \ \frac{2}{b}(\widetilde e_1 + \cdots + \widetilde e_{bk}), \,  \ldots \, , \, \frac{b-1}{b}(\widetilde e_1 + \cdots + \widetilde e_{bk}), 
\]
which have final coordinates $k, \, 2k, \ldots, (b-1)k$, respectively.
Therefore $B_{F,0}(t) = t^k + \cdots + t^{(b-1) k}$, and 
\[
h^*_P(t) = (1 + \cdots + t^d) + (1 + \cdots + t^r)(t^k + \cdots + t^{(b-1)k}).
\]
\end{example}

\begin{example}\label{dimension seven}
Suppose $d = 7$, and
\[
f = \frac{1}{7}(1,2,2,4,4,4,4).
\]
By Proposition~\ref{simplices}, $P = \conv\{ e_1, \ldots, e_7, -f \}$ is reflexive with respect to the lattice $N' = \Z^7 + \Z \cdot f$.  Let $F_1$, $F_2$, and $F_3$ be the faces of $\T$ with vertex sets $\{0, e_1, \ldots, e_7\}$, $\{0, e_1, e_2, e_3, -f\}$, and $\{0, e_1, -f\}$, respectively.  Then the lattice points in $\BBox (F_1, 0) \cap (N' \times \Z)$ are
\[
\begin{array}{l}
(1/7) ( \widetilde e_1 + 2 \widetilde e_2 + 2 \widetilde e_3 + 4 \widetilde e_4 + 4 \widetilde e_5 + 4 \widetilde e_6 + 4 \widetilde e_7), \\
 (1/7) ( \widetilde 2e_1 + 4 \widetilde e_2 + 4 \widetilde e_3 + \widetilde e_4 +  \widetilde e_5 + \widetilde e_6 + \widetilde e_7), \\
 (1/7) ( \widetilde 3e_1 + 6 \widetilde e_2 + 6 \widetilde e_3 + 5 \widetilde e_4 + 5 \widetilde e_5 + 5 \widetilde e_6 + 5 \widetilde e_7), \\
(1/7) ( 4 \widetilde e_1 + \widetilde e_2 + \widetilde e_3 + 2 \widetilde e_4 + 2 \widetilde e_5 + 2 \widetilde e_6 + 2 \widetilde e_7), \\
(1/7) ( 5 \widetilde e_1 + 3 \widetilde e_2 + 3 \widetilde e_3 + 6 \widetilde e_4 + 6 \widetilde e_5 + 6 \widetilde e_6 + 6 \widetilde e_7), \\ 
(1/7) ( 6 \widetilde e_1 + 5 \widetilde e_2 + 5 \widetilde e_3 + 3 \widetilde e_4 + 3 \widetilde e_5 + 3 \widetilde e_6 + 3 \widetilde e_7),
 \end{array}
 \]
 which have final coordinates 3, 2, 5, 2, 5, and 4, respectively.
So $B_{F_1,0}(t) = 2t^2 + t^3 + t^4 + 2t^5$.  Similarly, the lattice points in $\BBox (F_2, 0) \cap (N' \times \Z)$ are
\[
\begin{array}{l}
(1/4) (-3\widetilde f +  \widetilde e_1 + 2\widetilde e_2 + 2\widetilde e_3), \\
(1/4)(-\widetilde f + 3 \widetilde e_1 + 2\widetilde e_2 + 2\widetilde e_3),
\end{array}
\]
which both have final coordinate 2, so $B_{F_2,0}(t) = 2t^2$.  Finally, the unique lattice point in $\BBox (F_3, 0) \cap (N' \times \Z)$ is $(1/2)( -\widetilde f + e_1)$, so $B_{F_3, 0}(t)  = t$.  

For all other faces $F$, $\BBox(F,0)$ contains no lattice points.  It follows that
\[
h^*(P) = (1,2,6,5,5,6,2,1).
\]
\end{example}

\begin{example} \label{dimension eleven}
Suppose $d = 11$, and
\[
f = \frac{1}{11}(1,1,1,2,4,4,4,4,4,4,4).
\]
By Proposition~\ref{simplices}, $P$ is reflexive with respect to the lattice $N' = \Z^{11} + \Z \cdot f$.  Let $F_1$, $F_2$, and $F_3$ be the faces of $\T$ with vertex sets $\{0, e_1, \ldots, e_{11} \}$, $\{0, e_1, e_2, e_3, e_4, -f \}$, and $\{ 0, e_1, e_2, e_3, -f \}$, respectively.  Then, by computations similar to those in Example~\ref{dimension seven},
\begin{eqnarray*}
B_{F_1,0}(t) & =  & t^2 + 2t^3 + 2t^5 + 2t^6 + 2t^8 + t^9, \\
B_{F_2,0}(t) & = & t^2 + t^3, \\
B_{F_3,0}(t) & = & t^2,
\end{eqnarray*}
and $B_{F,0}(t) = 0$ for all other $F$.  It follows that
\[
h^*(P) = (1,1,4,6,4,6,6,4,6,4,1,1).
\]
\end{example}

\begin{proof}[Proof of Theorem~\ref{not unimodal}]
If $b \geq 3$ and $k > r+1$ in Example~\ref{reflexive simplex}, then $h^*(P)$ is not unimodal.  In particular, taking $b = 3$ and $r = 0, 1, \mbox{ or } 2$, produces reflexive simplices with nonunimodal $h^*$-vectors in all dimensions $d \geq 6$ except for $d = 7, 8, \mbox{ or } 11$.  For $d = 8$, one may take $b = 4$ and $k = 2$.  For dimensions $7$ and $11$, reflexive simplices with nonunimodal $h^*$-vectors are given by Examples~\ref{dimension seven} and \ref{dimension eleven}, respectively.
\end{proof}

\begin{proof}[Proof of Theorem~\ref{not macaulay}]
If $b = 2$ and $k \geq 2$ in Example~\ref{reflexive simplex}, then $h^*(P)$ is unimodal, but $h^*_1 - h^*_0 = 0$ and $h^*_k - h^*_{k-1} = 1$, so $g^*(P)$ is not a Macaulay vector.  Taking $b = k = 2$ produces such examples in all dimensions $d \geq 4$.
\end{proof}

One key ingredient in the proof of Theorem~\ref{big valleys} is the following special case of Braun's formula \cite{Braun06}, which gives the $h^*$-polynomial of the free sum of two reflexive polytopes.  Recall that if $Q$ and $Q'$ are polytopes in vector spaces $V$ and $V'$, respectively, each containing $0$ in its interior, then the free sum $Q \oplus Q'$ is the convex hull of $Q \times \{0 \}$ and $\{ 0 \} \times Q'$ in $V \times V'$.  

\begin{BF}
Let $Q$ and $Q'$ be reflexive polytopes.  Then
\[
h^*_{Q \oplus Q'}(t)  = h^*_Q(t) \cdot h^*_{Q'}(t).
\]
\end{BF}

\noindent Note that the dual of the free sum of two polytopes is the product of their respective duals, so the free sum of two reflexive polytopes is reflexive.

\begin{proof}[Proof of Theorem~\ref{big valleys}]
Nill has constructed a sequence of reflexive simplices $Q_j$ such that $\dim Q_j = j$ and the normalized volume
\[
\vol Q_j = h^*_0(Q_j) + \cdots + h^*_j(Q_j)
\]
grows doubly exponentially with $j$ \cite{Nill07}.  Therefore, there exists a reflexive simplex $Q$ such that $h^*_i(Q) \geq n$ for some $i$, and $\dim(Q) = O(\log \log n)$.  Let $Q'$ be the reflexive simplex constructed by taking $b = m + 1$, $k = \dim(Q) + 2$, and $d = bk$ in Example~\ref{reflexive simplex}.  Then $\dim(Q') = (m+1)(\dim(Q) + 2)$, so $\dim(Q \oplus Q') = O(m \log \log n)$, and
\[
h^*_{Q'}(t) = (1 + t + \cdots + t^d) + (t^k + t^{2k} + \cdots + t^{(b-1)k}).
\]

Let $P = Q \oplus Q'$.  Then $P$ is reflexive and, by Braun's Formula, $h^*_P(t) = h^*_Q(t) \cdot h^*_{Q'}(t)$.  It follows easily that 
\[
h^*_{k\ell + i}(P) = \vol(Q) + h_i^*(Q),
\]
for $1 \leq \ell \leq m+1$, and
\[
h^*_{k( \ell + 1) - 1}(P) = \vol(Q),
\]
for $1 \leq \ell \leq m$.  In particular, setting $i_\ell = k\ell + i$ and $j_\ell = k(\ell + 1) -1$, we have
$h^*_{i_\ell} - h^*_{j_\ell} \geq n$ and $h^*_{i_{\ell + 1}} - h^*_{j_\ell} \geq n$,
for $ 1 \leq \ell \leq m$, as required.
\end{proof}

\bibliography{math}

\providecommand{\bysame}{\leavevmode\hbox to3em{\hrulefill}\thinspace}
\providecommand{\MR}{\relax\ifhmode\unskip\space\fi MR }
% \MRhref is called by the amsart/book/proc definition of \MR.
\providecommand{\MRhref}[2]{%
  \href{http://www.ams.org/mathscinet-getitem?mr=#1}{#2}
}
\providecommand{\href}[2]{#2}
\begin{thebibliography}{10}

\bibitem{Athanasiadis04}
C.~Athanasiadis, \emph{{$h\sp \ast$}-vectors, {E}ulerian polynomials and stable
  polytopes of graphs}, Electron. J. Combin. \textbf{11} (2004/06), no.~2,
  Research Paper 6, 13 pp. (electronic).

\bibitem{Athanasiadis05}
\bysame, \emph{Ehrhart polynomials, simplicial polytopes, magic squares and a
  conjecture of {S}tanley}, J. Reine Angew. Math. \textbf{583} (2005),
  163--174.

\bibitem{BarvinokPommersheim99}
A.~Barvinok and J.~Pommersheim, \emph{An algorithmic theory of lattice points
  in polyhedra}, New perspectives in algebraic combinatorics (Berkeley, CA,
  1996--97), Math. Sci. Res. Inst. Publ., vol.~38, Cambridge Univ. Press,
  Cambridge, 1999, pp.~91--147.

\bibitem{BatyrevDais96}
V.~Batyrev and D.~Dais, \emph{Strong {M}c{K}ay correspondence, string-theoretic
  {H}odge numbers and mirror symmetry}, Topology \textbf{35} (1996), no.~4,
  901--929.

\bibitem{BeckRobins07}
M.~Beck and S.~Robins, \emph{Computing the continuous discretely: Integer point
  enumeration in polyhedra}, Undergraduate Texts in Mathematics,
  Springer-Verlag, 2007.

\bibitem{BeckSottile06}
M.~Beck and F.~Sottile, \emph{Irrational proofs for three theorems of
  {S}tanley}, European J. Combin. \textbf{28} (2006), no.~1, 403--409.

\bibitem{BetkeMcMullen85}
U.~Betke and P.~McMullen, \emph{Lattice points in lattice polytopes}, Monatsh.
  Math. \textbf{99} (1985), no.~4, 253--265.

\bibitem{Braun06}
B.~Braun, \emph{An {E}hrhart series formula for reflexive polytopes}, Electron.
  J. Combin. \textbf{13} (2006), no.~1, Note 15, 5 pp. (electronic).

\bibitem{BrunsRomer07}
W.~Bruns and T.~R{\"o}mer, \emph{{$h$}-vectors of {G}orenstein polytopes}, J.
  Combin. Theory Ser. A \textbf{114} (2007), no.~1, 65--76.

\bibitem{DeLoera05}
J.~De~Loera, \emph{The many aspects of counting lattice points in polytopes},
  Math. Semesterber. \textbf{52} (2005), no.~2, 175--195.

\bibitem{DeLoera04}
J.~De~Loera, R.~Hemmecke, J.~Tauzer, and R.~Yoshida, \emph{Effective lattice
  point counting in rational convex polytopes}, J. Symbolic Comput. \textbf{38}
  (2004), no.~4, 1273--1302.

\bibitem{Hibi90}
T.~Hibi, \emph{Some results on {E}hrhart polynomials of convex polytopes},
  Discrete Math. \textbf{83} (1990), no.~1, 119--121.

\bibitem{Hibi91}
\bysame, \emph{Ehrhart polynomials of convex polytopes, {$h$}-vectors of
  simplicial complexes, and nonsingular projective toric varieties}, Discrete
  and computational geometry (New Brunswick, NJ, 1989/1990), DIMACS Ser.
  Discrete Math. Theoret. Comput. Sci., vol.~6, Amer. Math. Soc., Providence,
  RI, 1991, pp.~165--177.

\bibitem{Hibi92}
\bysame, \emph{Algebraic combinatorics on convex polytopes}, Carslaw
  Publications, 1992.

\bibitem{Hibi94}
\bysame, \emph{A lower bound theorem for {E}hrhart polynomials of convex
  polytopes}, Adv. Math. \textbf{105} (1994), no.~2, 162--165.

\bibitem{MustataPayne05}
M.~Musta{\c{t}}{\v{a}} and S.~Payne, \emph{Ehrhart polynomials and stringy
  {B}etti numbers}, Math. Ann. \textbf{333} (2005), no.~4, 787--795.

\bibitem{Nill07}
B.~Nill, \emph{Volume and lattice points of reflexive simplices}, Discrete
  Comput. Geom. \textbf{37} (2007), no.~2, 301--320.

\bibitem{OhsugiHibi06}
H.~Ohsugi and T.~Hibi, \emph{Special simplices and {G}orenstein toric rings},
  J. Combin. Theory Ser. A \textbf{113} (2006), no.~4, 718--725.

\bibitem{Stanley80}
R.~Stanley, \emph{Decompositions of rational convex polytopes}, Ann. Discrete
  Math. \textbf{6} (1980), 333--342, Combinatorial mathematics, optimal designs
  and their applications (Proc. Sympos. Combin. Math. and Optimal Design,
  Colorado State Univ., Fort Collins, Colo., 1978).

\bibitem{Stanley93}
\bysame, \emph{A monotonicity property of {$h$}-vectors and {$h\sp
  *$}-vectors}, European J. Combin. \textbf{14} (1993), no.~3, 251--258.

\bibitem{Stienstra98}
J.~Stienstra, \emph{Resonant hypergeometric systems and mirror symmetry},
  Integrable systems and algebraic geometry (Kobe/Kyoto, 1997), World Sci.
  Publ., River Edge, NJ, 1998, pp.~412--452.

\end{thebibliography}
\bibliographystyle{amsplain}

\end{document}